\documentclass{amsart}
\usepackage{amsmath,amsthm,amssymb,IMjournal,graphicx}
\usepackage[linesnumbered]{algorithm2e}

\begin{document}
\newtheorem{The}{Theorem}[section]
\newtheorem{Def}{Definition}[section]
\newtheorem{Lem}{Lemma}[section]
\newtheorem{Pro}{Proposition}[section]
\newtheorem{Exa}{Example}[section]
\newtheorem{Alg}{Algorithm}
\input intmacros.sty

\numberwithin{equation}{section}

\title{Interval matrices with Monge property}

\author{\|Martin |\v{C}ern\'y|, Prague}

\rec {December 23, 2019}

\dedicatory{Cordially dedicated to Dr. Vladimir G. Deineko.}

\abstract 
   We generalize Monge property of real matrices for interval matrices. We define two classes of interval matrices with Monge property - in a strong and in a weak sense. We study fundamental properties of both classes. We show several different characterizations of the strong Monge property. For weak Monge property we give a polynomial characterization and several sufficient and necessary conditions. For both classes we study closure properties. We further propose a generalization of an algorithm by Deineko \& Filonenko which for a given matrix returns row and column permutations such that the permuted matrix is Monge if the permutations exist.
\endabstract

\keywords
   Monge matrix, interval matrix, interval analysis, linear programming
\endkeywords

\subjclass
65G99, 90C05, 15A57
\endsubjclass

\thanks
   The research has been supported by Czech Science Foundation Grant P403-18-04735S.
\endthanks

\section{Introduction}\label{sec1}
In 1781 a French mathematician Gaspard Monge observed a fundamental, but a very strong property while studying a variation of a transportation problem (see~\cite{Monge1781}). It was shown in the past century that the presence of Monge property (named in an honour to this great mathematician) simplifies many optimization problems. The famous NP-complete travelling salesman problem becomes solvable by a linear algorithm. Other optimization problems such as the assignment problem, the transportation problem or the lot-sizing problem can be solved significantly faster using algorithms based on Monge property. Since there is a geometrical interpretation of Monge property concerning distances, several applications in comuptational geometry are known. There are also further results in mathematical statistics, linguistics, bioinformatics, the graph theory or dynamic programming.\\
Interval analysis deals with an uncertainty or an inaccuracy in data. In almost every area of expertise people encounter a situation where they are limited by the precision of their data or their measuring devices. The problem becomes more severe when we use computers to compute abstract problems as a part of a mathematical proof. In these problems we cannot allow to neglect errors. In the interval analysis we envelope our data into intervals and then perform calculations on these intervals instead of the data itself. The methods of interval analysis ensure that the result is included in the resulting interval. In other typical problem for interval analysis we receive an interval of possible inputs and we want to find the range of all solutions.\\
Our work is the first study of an interval generalization of Monge property. For interval matrices we generalize the property in two natural ways - in a strong and in a weak sense. We show several characterizations of the interval matrices with the strong Monge property, few of them inspired by characterizations for real Monge matrices. We further state a polynomial characterization of matrices in the weak sense and study necessary and sufficient conditions. We also study closure properties of both classes of matrices. Finally, we present a permutation algorithm, which for a given general interval matrix decides if there exist row and column permutations such that the permuted matrix is Monge in the strong sense and returns the permutations if the answer is positive. 
\section{Preliminaries}\label{sec2}
\subsection{Interval analysis}
Before we start with an introduction to Monge matrices, we have to fix a notation and introduce basics of the interval analysis and interval arithmetics. For further information on the interval analysis see~\cite{AleMay2000},\cite{Hla2017da},\cite{GarAdm2016a}.\\
By $\mathbb{R}$ we denote the set of real numbers. We also denote by $\mathbb{IR}$ the set of all closed intervals over $\mathbb{R}$.
\begin{Def}(Interval matrix)\label{def1}
	An interval matrix $A \in \mathbb{IR}^{m\times n}$ is \[\A = \left[ \ul A, \ol A \right] = \left\{A \in \Rz^{m \times n}: \ul A \leq A \leq \ol A \right\}\] where $\ul A, \ol A$ are lower resp$.\ $upper bound matrices of $\A$.
\end{Def}
Similarly, we can define an interval vector as $\mathbf{v} = \left[ \ul v, \ol v \right] = \left\{v \in \Rz^{m}: \ul v \leq v \leq \ol v \right\}$.\\
Another way of defining an interval matrix is by using a so called \textit{center} $A^C = \frac{1}{2}(\ul A + \ol A)$ and a \textit{radius} $A^\Delta = \frac{1}{2} (\ol A - \ul A)$. Then an interval matrix can be rewritten as $\A = \left[A^C - A^\Delta,A^C + A^\Delta\right]$.
For two interval matrices $\M,\N \in \mathbb{IR}^{m\times n}$ we define intersection and union operations.
\begin{Def}(Interval matrix intersection)\label{def2}
	An interval matrix intersection $\M \cap \N$ is
	\[
	(\M \cap \N)_{ij}=
	\begin{cases}
	\left[l,u\right]& \text{if $l \leq u$},\\
	\emptyset& \text{if $l > u$},\\
	\end{cases}	
	\]
	where $l = \max\left\{\ul m_{ij},\ul n_{ij} \right\}$ and $u = \min\left\{\ol m_{ij},\ol n_{ij}\right\}$.
\end{Def}
\begin{Def}(Interval matrix union)\label{def3}
	For two interval matrices $\M,\N \in \mathbb{IR}^{m\times n}$ an interval matrix union is
	\(
	\M \cup \N = \left\{X \in \Rz^{m \times n} : X \in \M  \text{ or } X \in \N \right\}.
	\)
\end{Def}
Note that if $\M \cap \N = \emptyset$ then the interval matrix union is not an interval matrix. We deal with this by enveloping the set into an interval.
\begin{Def}(Envelope of interval matrix union)\label{def4}
	Let $\M \cup \N$ be an interval matrix union of two interval matrices $\M,\N \in \mathbb{IR}^{m\times n}$. An envelope of interval matrix union is $\square(\M \cup \N) = \left\{X \in \Rz^{m \times n} : \min\left\{\ul M, \ul N\right\} \leq X  \leq \max\left\{\ol M, \ol N\right\}\right\}.$
\end{Def}
\begin{Def}(Corner matrices)\label{def3.5}
	For \(\M \in \IR^{m \times n}\) an interval matrix, corner matrices $\downarrow M, \uparrow M$ are given by
\[
(\uparrow M)_{ij}=\left\{\begin{array}{l}
\ol m_{ij}\\
\ul m_{ij}
\end{array}\right\},
(\downarrow M)_{ij}=\left\{\begin{array}{l}
\ul m_{ij}\\
\ol m_{ij}
\end{array}\right\}
\text{ if } i+j \text{ is }
\left\{\begin{array}{l}
	\text{even}\\
	\text{odd}
\end{array}\right\}.
\]

\end{Def}
For a binary arithmetic operation $\circ \in \left\{+,-,\cdot,/\right\}$ defined on $\Rz$, we can introduce the corresponding interval operation as
\(
\a \circ \b = \left\{a \circ b : a \in \a, b \in \b \right\}.
\)
We can rewrite the definition into explicit formulae:
\begin{itemize}
	\item $\a + \b = \left[\ul a + \ul b, \ol a + \ol b\right]$,
	\item $\a - \b = \left[\ul a - \ol b, \ol a - \ul b\right]$,
	\item $\a \cdot \b = \left[\min\left\{\ul a\cdot\ul b, \ul a\cdot \ol b, \ol a \cdot \ul b, \ol a \cdot \ol b\right\},\max\left\{\ul a\cdot \ul b, \ul a \cdot \ol b, \ol a \cdot \ul b, \ol a \cdot \ol b\right\}\right]$,
	\item $\a / \b = \left[\min\left\{\ol a/\ol b, \ol a/\ul b,\ul a/\ol b,\ul a/\ul b\right\},\max\left\{\ol a/\ol b, \ol a/\ul b,\ul a/\ol b,\ul a/\ul b\right\}\right]$ if $0 \notin \b.$
\end{itemize}
Let us note that for the interval division there is a known generalization where $0 \in \b$.
\subsection{Real matrices with Monge property}
All the results from this subsection can be found in a survey by Burkard (see~\cite{persp}).
\begin{Def}(Monge matrix)\label{def5}
	Let $M \in \Rz^{m \times n}$. The matrix $M$ is Monge if for all $i,j,k,\ell :$ $1 \leq i < k \leq m$, $1 \leq j < \ell \leq n$ it holds \[m_{ij} + m_{k\ell} \leq m_{i\ell} + m_{kj}.\]
\end{Def}
Since Hoffman rediscovered Monge property in 1961, several equivalent characterizations have been shown. We merge some of the characterizations into a theorem, but first, we define a notion of \textit{submodular} functions.\
\begin{Def}(Submodular function)\label{def6}
	Let $\Lambda = (I,\wedge,\vee)$ be a distributive lattice where $I = \left\{1,...,m\right\}\times\left\{1,...,n\right\}$ and join ($\wedge$)  and meet ($\vee$) operations are defined for $x = (x_1,x_2) ,y = (y_1,y_2)$ as follows:
	\begin{itemize}
		\item $(x_1,x_2) \wedge (y_1,y_2)  = (\min\left\{x_1,y_1\right\},\min\left\{x_2,y_2\right\})$,
		\item $(x_1,x_2) \vee (y_1,y_2)  = (\max\left\{x_1,y_1\right\},\max\left\{x_2,y_2\right\})$.
	\end{itemize}
	Function $f: I \rightarrow \Rz$ is said to be submodular on $\Lambda$ if for all $x,y \in I$
	\[f(x \vee y) + f(x \wedge y) \leq f(x) + f(y).\]
\end{Def}
\begin{The}\label{thm1}
	Let $M \in \Rz^{m \times n}$, then the following are equivalent:
	\begin{enumerate}
		\item $M$ is Monge matrix,
		\item $m_{ij} + m_{k\ell} \leq m_{i\ell} + m_{kj}$ for all $1 \leq i < k \leq m, 1 \leq j < \ell \leq n$,
		\item $m_{ij} + m_{i+1,j+1} \leq m_{i,j+1} + m_{i+1,j}$ for all $1 \leq i  < m, 1 \leq j < n$,
		\item A function $f: I \rightarrow \Rz$ defined by $f(i,j) = m_{ij}$ is submodular on $\Lambda$ where $\Lambda = (I,\wedge,\vee)$ is a distributive lattice.
	\end{enumerate}
\end{The}
We further present a list of operations under which Monge matrices are closed.
\begin{The}\label{thm2}
	Let $M,N \in \mathbb{R}^{m \times n}$ be Monge. Then the following holds:
	\begin{enumerate}
		\item $M^T$ is Monge,
		\item $\alpha M$ is Monge for $\alpha \geq 0$,
		\item $M + N$ is Monge,
		\item for any $u \in \Rz^{m}$, $v \in \Rz^{n}$, matrix $C \in \mathbb{R}^{m\times n}$ defined by $c_{ij} = m_{ij} + u_i + v_j$ is Monge.
	\end{enumerate}
\end{The}
The second and the third result in Theorem~\ref{thm2} imply that the set of nonnegative Monge matrices forms a convex polyhedral cone. This cone can be described by 4 types of 0-1 matrices corresponding to the extremal rays. Let $H^{i}$ denote a 0-1 matrix where $i$th row contains all ones while the other entries are zeros and $V^{j}$ a 0-1 matrix with $j$th columns set to ones and the rest to zeros. Further, let $L^{rs}$ be a 0-1 matrix where for $l_{ij}^{rs} = 1$ for $i=r,\dots,m$ and $j=1,\dots,s$. Otherwise $l_{ij}^{rs}=0$. Similarly let $R^{pq}$ be a 0-1 matrix with $r_{ij}^{pq}=1$ for $i=1,\dots,p$ and $j=q,\dots,n$, otherwise $r_{ij}^{pq}=0$. Any Monge matrix can be represented by a nonnegative combination of matrices $H^{i},V^{j},L^{rs}$ and $R^{pq}$. 
\begin{The}\label{thm3}
	Let $M \in \mathbb{R}^{m\times n}$ be Monge matrix, then there are coefficients $\kappa_{i},\lambda_{j},\mu_{rs}$ and $\nu_{pq}$ such that
	\[M = \sum_{i=1}^{m}\kappa_{i}H^{i} + \sum_{j=1}^{n}\lambda_jV^{j}+\sum_{r=2}^m\sum_{s=1}^{n-1}\mu_{rs}L^{rs}+ \sum_{p=1}^{m-1}\sum_{q=2}^n\nu_{pq}R^{pq}.\]
	The matrices $H^{p}$ with $p=1,\dots,m$, $V^{q}$ with $q=1,\dots,n$, $L^{rs}$ with $r=2,\dots,m,$ $s=1,\dots,n-1$ and $R^{pq}$ with $p=1,\dots,m-1,q=2,\dots,n$ generate extreme rays of the cone of nonnegative Monge matrices.
\end{The}
\section{Interval matrices with strong Monge property}
In this section we introduce interval matrices with the strong Monge property. We present a list of 5 equivalent characterizations, most of them similar to those in Theorem~\ref{thm1}.
\begin{Def}(Strong Monge property)\label{def7}
	An interval matrix $\M \in \mathbb{IR}^{m\times n}$ has the strong Monge property if every $M \in \M$ is Monge. We denote by $\mathbb{ISM}$ the set of interval matrices with strong Monge property.
\end{Def}
Before we state the equivalent characterizations we first need to define a generalization of submodular functions.
\begin{Def}(Interval submodular functions)\label{def8}
	Let $\Lambda = (I,\wedge,\vee)$ be a distributive lattice where $I = \left\{1,...,m\right\}\times\left\{1,...,n\right\}$ with join $(\wedge)$ and meet $(\vee)$ operations. The operations are defined for $x = (x_1,x_2), y = (y_1,y_2)$ as follows:
	\begin{itemize}
		\item $(x_1,x_2) \wedge (y_1,y_2) = (\min\left\{x_1,y_1\right\},\min\left\{x_2,y_2\right\})$,
		\item $(x_1,x_2) \vee (y_1,y_2) = (\max\left\{x_1,y_1\right\},\max\left\{x_2,y_2\right\})$.
	\end{itemize}
	A function $\f : I \rightarrow \IR$ is \textit{submodular on lattice $\Lambda$} if $\ol f(x\vee y) + \ol f(x \wedge y) \leq \ul f(x) + \ul f(y)$ for all $x,y \in I$.
\end{Def}

\begin{The}(Characterization of strong Monge property)\label{thm4}
	Let $\M \in \IR^{m\times n}$ be an interval matrix. Then the following are equivalent:
	\begin{enumerate}
		\item $\M \in \ISM$,
		\item $\overline{m}_{ij} + \overline{m}_{k\ell} \leq \underline{m}_{i\ell} + \underline{m}_{kj}$ for all $1 \leq i < k \leq m, 1 \leq j < \ell \leq n$,
		\item $\overline{m}_{ij} + \overline{m}_{i+1,j+1} \leq \underline{m}_{i,j+1} + \underline{m}_{i+1,j}$ for all $1 \leq i  < m, 1 \leq j < n$,
		\item Corner matrices $\downarrow M$ and $\uparrow M$ are Monge,
		\item A function $\f : I \rightarrow \IR$ defined by $\f(i,j) = \m_{ij}$ is submodular on $\Lambda$ where $\Lambda = (I,\wedge,\vee)$ is a distributive lattice.
	\end{enumerate}
\end{The}
\proof
$(1) \leftrightarrow (2) \leftrightarrow (3)$: can be easily derived using Theorem~\ref{thm1} and Definition~\ref{def7}.\\
$(3) \leftrightarrow (4)$: can be derived using Definition~\ref{def7} and Definition~\ref{def3.5}.\\
$(3) \leftrightarrow (5)$:
Let $\M \in \mathbb{IR}^{m\times n}$ such that $(3)$ holds. Let further $x = \left(i,j+1\right) \in I$ and $y = \left(i+1,j\right) \in I$. Then 
\[
\ol f(x \wedge y) = \ol f(\left(i,j+1\right) \wedge \left(i+1,j\right)) = \ol f((i,j)) = \ol m_{ij}
\]
and
\[
\ol f(x \vee y) = \ol f(\left(i,j+1\right) \vee \left(i+1,j\right)) = \ol f((i+1,j+1)) = \ol m_{i+1,j+1}.
\]
Therefore
\[
\ol f(x \vee y) + \ol f(x \wedge y) = \ol m_{i+1,j+1} + \ol m_{ij} \leq \ul m_{i,j+1} + \ul m_{i+1,j} = \ul f(x) + \ul f(y).
\]
Since the inequality holds for any $i,j$, $\f$ is submodular on $\Lambda$.\\
Let us now suppose that the function $\f$ is submodular on the lattice $\Lambda$. Then the condition 
\[
\ol f((i+1,j) \wedge (i,j+1)) + \ol f((i+1,j) \vee (i,j+1)) \leq  \ul f((i+1,j)) + \ul f((i,j+1)).
\]
corresponds to
\[
\ol m_{ij} + \ol m_{i+1,j+1} \leq \ul m_{i+1,j} +  \ul m_{i,j+1}
\]
for every $i,j$, thus $\M \in \ISM$.

\endproof
Let us remark that the result in Theorem~\ref{thm3} does not seem to be easily generalizable to the interval case. Trying to find an \textit{interval decomposition} by taking one possible decomposition for each $M \in \M$ and making an interval envelope of all possible coefficients $\kappa_i,\lambda_j,\mu_{rs},\nu_{pq}$ leads to an overestimation in general as shown in the example below.
\begin{Exa}
	Let $\M \in \ISM$ such that\[
	\M=
	\begin{pmatrix}
	[0,5] & 5\\
	[0,8] & 0 \\
	\end{pmatrix}.\]
	If the decomposition is to equal $\M$ then it must be in the form of
	\[\C(\M)=
	[0,3]
	\begin{pmatrix}
	0 & 0\\
	1 & 0\\
	\end{pmatrix}
	+[0,5]
	\begin{pmatrix}
	1 & 0\\
	1 & 0\\
	\end{pmatrix}
	+ 5
	\begin{pmatrix}
	0 & 1\\
	0 & 0\\
	\end{pmatrix}.\]
	But for matrix 
	\[M = \begin{pmatrix}
	1 & 5\\
	6 & 0\\
	\end{pmatrix}
	\]
	we see that there is no possible decomposition of $M$ between the coefficients of $\C(\M)$.
\end{Exa}
The described overestimating decomposition can be computed by the linear programming but since we do not need it further in the text, we omit the construction.
\section{Interval matrices with the weak Monge property}
In this section we introduce the interval matrices with the weak Monge property. We offer a polynomial characterization and several necessary and sufficient conditions.
\begin{Def}(Weak Monge property)
	An interval matrix $\M \in \mathbb{IR}^{m\times n}$ has the weak Monge property if there is Monge matrix $M \in \M$. We denote by $\mathbb{IWM}$ the set of interval matrices with the weak Monge property.
\end{Def}
We start off by showing that matrices with the weak Monge property are polynomially recognizable by a reduction to a special linear program.
\begin{The}\label{thm5}
	Let $\M \in \IR^{m \times n}$ and let $LP(\M)$ be a linear program defined as
	
	\begin{equation*}
	\begin{array}{llrr}
	\text{minimize}	 & const.\\	
	\text{subject to} & m_{ij} + m_{i+1,j+1} - m_{i,j+1} - m_{i+1,j} \leq 0, &(1)\\
	& m_{k\ell } \leq \ol m_{k\ell },&(2)\\
	& - m_{k\ell } \leq  - \ul m_{k\ell },& (3)\\
	\\
	\text{where} &1 \leq i < m, \text{ } 1 \leq j < n,\\
	&1 \leq k \leq m, \text{ } 1 \leq \ell \leq n.\\
	\end{array}
	\end{equation*}
	Then the matrix $\M \in \IWM$ iff $LP(\M)$ has a feasible solution.
\end{The}
\proof
	A feasible solution of $LM(\M)$ corresponds to a matrix $M$. Monge property of the matrix is guaranteed by $(1)$ and by $(2)$ and $(3)$ every entry $\ul m_{ij} \leq m_{ij} \leq \ol m_{ij}$ is from $\m_{ij}$.
	Thus every feasible solution of $LP(\M)$ is Monge matrix $M \in \M$ and therefore $\M \in \IWM$. If the linear program is not feasible, $\M$ does not have the weak Monge property.
\endproof
	Theorem~\ref{thm5} is important because we know that the recognition problem of matrices with the weak Monge property is solvable in polynomial time~\cite{GreHoo95}. For $\IWM$ we did not find any other polynomial characterization. Let us note that all of the characterizations of real Monge matrices can be restated for $\IWM$, although none of them can be used without any further modification to construct an efficient polynomial recognition algorithm.
\subsection{Necessary conditions}
Although we know the recognition problem of $\IWM$ is polynomial, the only characterization we found was by linear programming which is categorized as one of the hardest problems in the hierarchy of polynomially solvable problems (see \cite{GreHoo95}). Therefore we investigated necessary and sufficient conditions of~$\IWM$.\\
The first necessary condition employs so called \textit{residual} matrices.
\begin{Def}\label{def10}
Let $\M \in \IR^{m\times n}$ be an interval matrix. Then an interval residual matrix $\M^R \in \IR^{(m-1)\times (n-1)}$ is defined as
\[
\m^R_{ij} =  \left[\ul m_{i+1,j} + \ul m_{i,j+1} - \ol m_{ij} - \ol m_{i+1,j+1}, \ol m_{i+1,j} + \ol m_{i,j+1} - \ul m_{ij} - \ul m_{i+1,j+1} \right].
\]
\end{Def}
The residual matrices carry an information about the tightness of inequalities from the definition of Monge property.
\begin{Pro}\label{thm6}
	Let $\M \in \IWM$ and $\M^{R}$ be its residual matrix. Then there exists
	$M^R \in \M^R$ such that $M^R$ is nonnegative.
\end{Pro}
\proof
Let $\M \in \IWM$ and $M \in \M$ such that $M$ is Monge. Because of Monge property of $M$ we have $m^R_{ij} = m_{i+1,j} + m_{i,j+1} - m_{ij} - m_{i+1,j+1} \geq 0$ for all $i,j$. If we take matrix $(M^R)_{ij} = m^R_{ij}$ it clearly holds $M^R \in \M^R$.
\endproof
Another necessary condition considers a presence of a special Monge matrix in the correspoing interval matrix with the weak Monge property.
\begin{Pro}\label{thm7}
	Let $\M \in \IWM^{m \times n}$. Then there exists $M \in \M$ such that $M$ is Monge and the number of entries $m_{ij} = \ol m_{ij}$ is at least $\max\left\{m,n\right\}$.
\end{Pro}
\proof
Let $M \in \M$ be Monge matrix. By Theorem~\ref{thm3} we can rewrite $M$ as 
\[
M = \sum_{i=1}^{m}\kappa_{i}H^{i} + \sum_{j=1}^{n}\lambda_jV^{j}+\sum_{r=2}^m\sum_{s=1}^{n-1}\mu_{rs}L^{rs}+ \sum_{p=1}^{m-1}\sum_{q=2}^n\nu_{pq}R^{pq}.
\]
	Let us take $M$ such that the number of entries $m_{ij} = \ol m_{ij}$ in $M$ is the highest possible and still lower than $\max\left\{m,n\right\}$. Let us suppose that $m > n$. It means that there is a row $k$ in $M$ where $m_{kj}\neq \ol m_{kj}$ for every column $j$. We take $\mu = \min\limits_{j}\left\{\ol m_{kj} - m_{kj} \right\}$ and add $\mu H^{i}$ to $M$. The matrix $M + \mu H^{i}$ is also Monge, belongs to $\M$ and the number of upper bounds of intervals in $M + \mu H^{i}$ is higher than in $M$.\\
	For $n > m$ we employ the matrices of type $V^{j}$ and the rest of the argument is similar.
\endproof
	To show that the bound in Proposition~\ref{thm7} can be achieved we give the following example.
\begin{Exa}
			Let $\M \in \IR^{4 \times 4}:$
	\[\M = \begin{pmatrix}
	\left[3,1000 \right] & \left[10,120 \right] & \left[17,20 \right] & \left[0,24 \right] \\
	\left[2,20 \right] & \left[7,9 \right] & \left[0,12 \right] & \left[17,85 \right] \\
	\left[2,5 \right] & \left[0,6 \right] & \left[10,14 \right] & \left[14,100 \right]\\
	\left[0,1 \right] & \left[3,6 \right] & \left[5,21 \right] & \left[7,1000 \right]\\
	\end{pmatrix}.
	\]
	Matrix $M \in \M$ such that
	\[M = \begin{pmatrix}
		3 & 10 & 17 & 24 \\
	2 & 7 & 12 & 17 \\
	2 & 6 & 10 & 14\\
	1 & 3 & 5 & 7 \\
	\end{pmatrix}.
	\]
	is Monge, therefore $\M \in \IWM$. Moreover, on the diagonal from the lower left corner to the upper right corner the values are upper bounds of the corresponding intervals. It is easy to check that for any Monge matrix $N \in \M$, no other entry can be an upper bound of $\M$ since it would violate at least one of neighbouring conditions of Monge property.
\end{Exa}		

\subsection{Sufficient conditions of matrices with the weak Monge property}
The first two sufficient conditions use the decomposition into extremal rays of convex cone (see Theorem~\ref{thm3}).
\begin{Pro}\label{thm8}
	Let $\M \in \IR^{m \times n}$. If it holds for every row $i$ that $\bigcap\limits_{j} \m_{ij} \neq \emptyset$ or for every column $j$ that $\bigcap\limits_{i} \m_{ij} \neq \emptyset$, then $\M \in \IWM$.
\end{Pro}
\proof
	Let us suppose that for every row $i$ it holds that $\bigcap\limits_{j} \m_{ij} = \left[\ul \alpha_i, \ol \alpha_i\right]$. Then a matrix
	\[
	M = \alpha_1H^{1} + \alpha_2H^{2} + \dots + \alpha_nH^{n}
	\]
	where $\alpha_i \in \left[\ul \alpha_i, \ol \alpha_i \right]$ is Monge matrix by Theorem~\ref{thm3}. Since $M \in \M$, we conclude that $\M \in \IWM$.
	For nonempty intersections of columns the argument is similar.
\endproof

\begin{figure}[h]
		\includegraphics{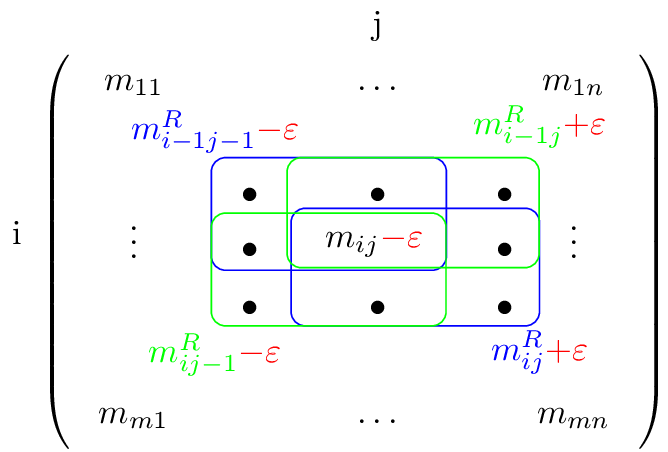}
		\caption{Subtracting $\varepsilon$ from $m_{ij}^C$ and its effect on entries of $M^R$.}\label{fig1}
	\end{figure}

\begin{figure}[h]
		\includegraphics{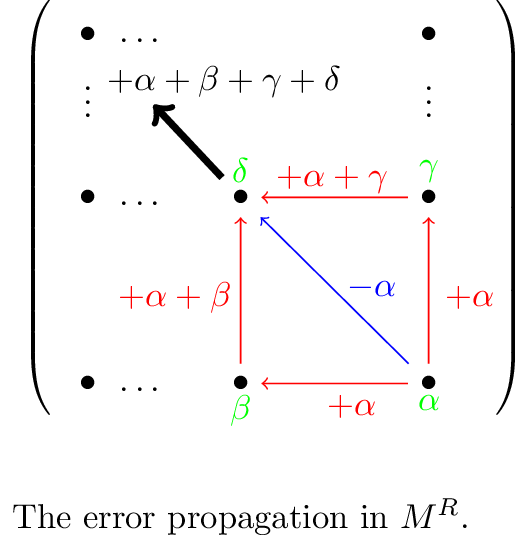}
		\caption{The residual propagation in $M^R$.}\label{fig2}
	\end{figure}

\begin{figure}[h]
		\includegraphics{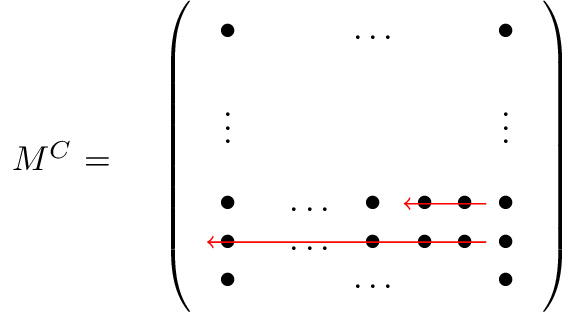}
		\caption{The order of changing values in $M^C$ to zero the entries of $M^R$.}\label{fig3}
	\end{figure}

\begin{Pro}\label{thm9}
	Let $\M \in \IR^{m \times n}$. If it holds for all indices $i,j$ that $m_{ij}^{\Delta}\geq \lvert m_{ij}^{C}$ then $\M \in \IWM$.
\end{Pro}
\proof
	The condition $m_{ij}^{\Delta}\geq \lvert m_{ij}^{C}$ is equivalent with $0 \in \m_{ij}$. Thus Monge matrix $0_{m \times n} \in \M$.
\endproof
	Another class of sufficient conditions of matrices with the weak Monge property is based on an idea that in a space of real matrices we start with $M^C$ and use an easy procedure to move in steps from $M^C$ until we reach Monge matrix of a special form. Depending on the direction and distance of each step we can compute how far we have to move from $M^C$ in each interval entry to achieve the matrix. By this, we can get a sufficient condition dependent on the width of intervals. To determine the necessary width of intervals we employ residual matrices.
\begin{The}\label{thm10}
	Let $\M \in \IR^{m \times n}$ and let $M^R \in \Rz^{(m-1)\times (n-1)}$ be the residual matrix of $M^C$ meaning $(M^R)_{ij} = m_{i+1,j} + m_{i,j+1} - m_{ij} - m_{i+1,j+1}$. If for all indices $i,j$ of $\M$ it holds that $m^\Delta_{ij} \geq \lvert\sum\limits_{k=i}^{m-1}\sum\limits_{\ell=j}^{n-1}m^R_{k\ell}\rvert$ then $\M \in \IWM$.
\end{The}
\proof
	Let $M^C \in \Rz^{m \times n}$ and let $M^R \in \Rz^{(m-1)\times (n-1)}$ be its residual matrix. In general, the residual matrix $M^R$ will not be nonnegative. Our goal is to set the entries of $M^R$ to zero by changing the entries of $M^C$. We set the entries to zero one by one using a specific elimination order. We see that by subtracting $\varepsilon$ from $m^C_{ij}$ the value of $m^R_{ij}$ increases by $\varepsilon$. By this operation, entries $m^R_{i-1,j-1},m^R_{i-1,j}$ and $m^R_{i,j-1}$ are affected as well (see Figure~\ref{fig1}).
	We start from the bottom-right corner of $M^R$ and add the value of $m^R_{m-1,n-1}$ to $m^C_{m-1,n-1}$. This sets the residuum $m^R_{m-1,n-1}$ to zero and propagates its value into the three neighbouring entries (see Figure~\ref{fig2}).
	In next step, we eliminate the residuum of the element $m^R_{m-1,n-2}$ and continue in the decreasing order of columns until we arrive at the beginning of the row, then proceed with the row above in the same manner (see Figure~\ref{fig3}).\\
	By each step we eliminate one residuum and more importantly, no residuum already set to zero is affected further in the process.\\
	Not only this elimination order yields $0^{(m-1) \times (n-1)}$ residual matrix (therefore a corresponding Monge matrix) but it is also easy to describe the propagation of residual values in $M^R$. Setting to zero the residuum $m^R_{ij} = \alpha$ adds $\alpha$ to $m^R_{i-1,j}$ and $m^R_{i,j-1}$ and subtracts it from $m^R_{i-1,j-1}$ (as shown in Figure~\ref{fig2}). Now if the intervals of $\M$ are large enough we can move from the central matrix $M^C$ far enough to eliminate the residua. It is now easy to compute by induction the necessary condition for each interval of $\M$.\\

\begin{figure}[h]
	\includegraphics{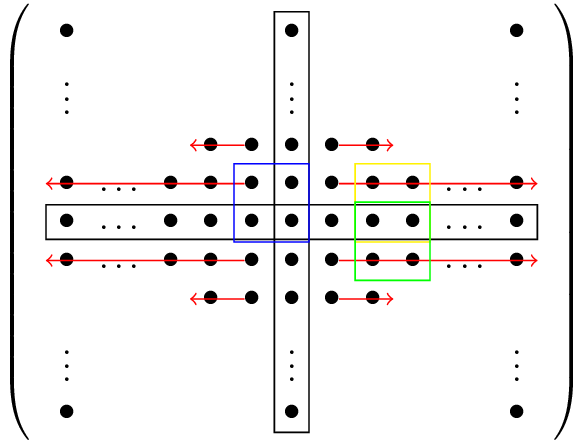}
	\caption{The residual propagation does not interfere between the blocks.}
	\label{fig4}
\end{figure}

	For the base step, from the way of propagation (ilustrated by Figure~\ref{fig2}) it is clear that it must hold that
	\begin{itemize}
		\item $m^\Delta_{m-1,n-1} \geq \lvert m^R_{m-1,n-1}\rvert$
		\item $m^\Delta_{m-1,n-2} \geq \lvert m^R_{m-1,n-1} + m^R_{m-1,n-2}\rvert$,
		\item $m^\Delta_{m-2,n-1} \geq \lvert m^R_{m-1,n-1} + m^R_{m-2,n-1}\rvert$,
		\item $m^\Delta_{m-2,n-2} \geq \lvert m^R_{m-2,n-1} + m^R_{m-1,n-2} +2m^R_{m-1,n-1} - m^R_{m-1,n-1}\rvert$, therefore\\ $m^\Delta_{m-2,n-2} \geq \lvert m^R_{m-2,n-1} + m^R_{m-1,n-2} + m^R_{m-1,n-1}\rvert$.
	\end{itemize}
	For inductional step let us suppose the residuum $m^R_{ij}$. It must hold that
	\[m^\Delta_{ij} \geq \left\lvert m^R_{ij} + m^R_{i+1,j} +  m^R_{i,j+1} - m^R_{i+1,j+1}\right\rvert.\]
	By induction we know that the residues are equal to
	\[
	m^\Delta_{ij} \geq\left\lvert m^R_{ij} + \sum\limits_{k=i+1}^{m-1}\sum\limits_{\ell=j}^{n-1}m^R_{k\ell} +  \sum\limits_{k=i}^{m-1}\sum\limits_{\ell=j+1}^{n-1}m^R_{k\ell} - \sum\limits_{k=i+1}^{m-1}\sum\limits_{\ell=j+1}^{n-1}m^R_{k\ell}\right\rvert
	\]
	which is equal to the form stated in the theorem.
\endproof
Let us note that the condition we just showed can be checked in $O(mn)$ time using dynamic programming.\\
The sufficient condition shown in the previous theorem is one of many modifications of the same condition depending on the order we choose to zero the values in $M^R$. The advantage of this one-diagonal order is that it is easy to compute the width of intervals.
We present one more condition from this class. The previous condition works well when the sum $\lvert\sum\limits_{k=i}^{m-1}\sum\limits_{l=j}^{n-1}m^R_{ij}\rvert \sim 0$ or is at least small for every $i,j$. If the errors are of the same sign, however, the sum has a tendency to grow a lot. This is because we propagate the error only in one direction.\\
We can choose a point in the matrix and propagate the error in four different (diagonal) directions.\\
\begin{The}\label{thm11}
	Let $\M \in \IR^{m \times n}$ and let $M^R \in \Rz^{(m-1)\times (n-1)}$ be the residual matrix of $M^C$. If there are indices $i,j$ of $\M$ such that
	\begin{itemize}
		\item $m_{rs}^{\Delta} \geq \lvert \sum\limits_{k=r}^{i-1}\sum\limits_{\ell=s}^{j-1} m_{k\ell}^R\rvert \text{ for every } (r < i) \wedge (s < j),$
		\item $m_{rs}^{\Delta} \geq \lvert \sum\limits_{k=r}^{i-1}\sum\limits_{\ell=j}^{s-1} m_{k\ell}^R\rvert \text{ for every } (r < i) \wedge (s > j),$
		\item $m_{rs}^{\Delta} \geq \lvert \sum\limits_{k=i}^{r}\sum\limits_{\ell=s}^{j-1} m_{k\ell}^R\rvert \text{ for every } (r > i) \wedge (s < j),$
		\item $m_{rs}^{\Delta} \geq \vert \sum\limits_{k=i}^{r}\sum\limits_{\ell=j}^{s} m_{k\ell}^R\rvert \text{ for every } (r > i) \wedge (s > j),$
	\end{itemize} 
	then $\M \in \IWM$.
\end{The}
\proof
Let $i,j$ be indices of $M^C$. Then we can take $m^R_{i-1,j-1},m^R_{i-1,j+1},m^R_{i+1,j-1}$ and $m^R_{i+1,j+1}$ as starting points for residual elimination described in Theorem~\ref{thm10}. We can see in Figure~\ref{fig4} that the residua are not propagated between the blocks of $M^R$. The inequalities follow from Theorem~\ref{thm10}.
\endproof

\section{Closure properties of interval matrices with Monge property}
In this section we briefly introduce closure properties of both classes of interval matrices as well as those interconnecting them.
\subsection{Closure properties of matrices with the strong Monge property}
As mentioned in preliminaries the set of nonnegative real Monge matrices forms a convex cone meaning the matrices are closed under linear combinations with nonnegative coeficients. The fact that matrices with the strong Monge property are convex subsets of the set of real Monge matrices promises similar results for $\mathbb{ISM}$.
\begin{Pro}
	Let $\M,\N \in \mathbb{ISM}$ and let $\alpha \in\mathbb{R}^{+}_0$. Then also $\alpha \M \in \ISM$ and $\M + \N \in \ISM$.
\end{Pro}
\proof
Straightforward from Theorem~\ref{thm4}.(3).
\endproof
When it comes to multiplication by interval \boldmath$\alpha$ \unboldmath $\in \IR^{+}_0$, interval matrices with strong Monge property are closed only under certain restriction dependent on the lower bound of \boldmath$\alpha$ \unboldmath and its radius.
\begin{The}
	Let $\M \in \ISM^{+}_0$ and let \(\boldsymbol\alpha\in \IR_{0}^{+}\). Then \boldmath$\alpha$\unboldmath $\M \in \ISM^{+}_0$ iff
	\[
	\frac{\alpha^{\Delta}}{\alpha^C} \leq \phi \mbox{ where } \phi =\min\limits_{i,j}\left(\frac{\ul m_{i,j+1} + \ul m_{i+1,j} - \ol m_{ij} - \ol m_{i+1,j+1}}{\ul m_{i,j+1} + \ul m_{i+1,j} + \ol m_{ij} + \ol m_{i+1,j+1}}\right).
	\] 
\end{The}

\proof
For all indices $i,j$ it must hold that\\
\[
\boldsymbol\alpha\ol m_{ij} + \boldsymbol\alpha\ol m_{i+1,j+1} \leq \boldsymbol\alpha\ul m_{i,j+1} +\boldsymbol\alpha\ul m_{i+1,j}.
\] 
It holds for all $\alpha \in \boldsymbol\alpha$ that
\[ 
\alpha\ol m_{ij} + \alpha\ol m_{i+1,j+1}\leq
\ol \alpha \mbox{ } \ol m_{ij} + \ol \alpha \mbox{ }  \ol m_{i+1,j+1} \leq \ul \alpha \mbox{ } \ul m_{i,j+1} + \ul \alpha \mbox{ }  \ul m_{i+1,j}
\leq \alpha\ol m_{i,j+1} + \alpha\ol m_{i+1,j}.
\]
We achieve the tightest inequality for
\[
\ol \alpha \mbox{ } \ol m_{ij} + \ol \alpha \mbox{ }  \ol m_{i+1,j+1} \leq \ul \alpha \mbox{ } \ul m_{i,j+1} + \ul \alpha \mbox{ }  \ul m_{i+1,j}.
\]
Adjusting the inequality, we get
\[
\overline{\alpha} \leq \underline{\alpha}\left(\frac{\underline{m}_{i,j+1} + \underline{m}_{i+1,j}}{\overline{m}_{ij} + \overline{m}_{i+1,j+1}}\right).
\]
Substituting $\ol \alpha$ for $\alpha^C + \alpha^{\Delta}$, $\ul \alpha$ for $\alpha^C - \alpha^\Delta$ and adjusting again the inequality we get the formula
\begin{equation}\label{eq02:1}
\frac{\alpha^\Delta}{\alpha^C} \leq \left(\frac{\ul m_{i,j+1} + \ul m_{i+1,j} - \ol m_{ij} - \ol m_{i+1,j+1}}{\ul m_{i,j+1} + \ul m_{i+1,j} + \ol m_{ij} + \ol m_{i+1,j+1}}\right).
\end{equation}
It is now clear that the inequality~\ref{eq02:1} holds for all $i,j$ iff it holds for minimum over all indices.
\endproof
Finally, we state two observations. The first one is about matrix transposition and the second one about matrix products.
\begin{Pro}
	For a matrix $\M \in \ISM$ the transposition $\M^T \in \ISM$.
\end{Pro}
\proof
	Straightforward from the definition of $\ISM$.
\endproof
\begin{Exa}
	Let us consider matrices
	\[A = \begin{pmatrix}
	5 & 5 \\
	0.1 & 0.1
	\end{pmatrix}
	B = \begin{pmatrix}
	5 & 0.1 \\
	6 & 0.1
	\end{pmatrix}.\]
	The matrix $A\odot B \notin \ISM$ for $\odot$ representing the Standard, the Hadamard and the Kronecker (tensor) matrix product.
\end{Exa}
\subsection{Closure properties of matrices with the weak Monge property}
	We investigated closure properties of several operations on $\IWM$. Most of the results are easy to prove, therefore we state them in one theorem.
	
	\begin{The}
		Let $\P \in \IR^{m \times n}, \M,\N \in \IWM^{m \times n}$, $\alpha \in \Rz^{+}_0$ and $\boldsymbol{\alpha} \in \IR^{+}_0.$ Then the following holds.
		\begin{enumerate}
			\item $\M + \N \in \IWM$,
			\item $\M + \P \in \IWM$ iff $\overline{\M^R} +  \overline{\P^R} \geq 0$,
			\item $\square (\M \cup \P) \in \IWM$,
			\item $\alpha\M \in \IWM$,
			\item $\boldsymbol{\alpha}\M \in \IWM$.
		\end{enumerate}
	\end{The}	
	\proof
		All the results are easy to prove from the definition of $\IWM$.
	\endproof
\subsection{Closure properties interconnecting both classes}	
	\begin{The}
		Let $\M \in \ISM^{m \times n},\N \in \IWM^{m \times n}$, $\alpha \in \Rz^{+}_0$ and $\boldsymbol{\alpha} \in \IR^{+}_0.$ Then the following holds.
		\begin{enumerate}
			\item $\M + \N \in \IWM$,
			\item $\forall i,j$ it holds that $\m_{ij} \cap \n_{ij} \neq \emptyset \rightarrow \M \cap \N \in \IWM$,
			\item $\square (\M \cup \N) \in \IWM$.
		\end{enumerate}
	\end{The}
	\proof
		All the results are easy to prove from the definition of $\IWM$ and $\ISM$.
	\endproof
	
\section{Permutation algorithm for Monge permutable matrices}
In many optimization problems (e.g. the travelling salesman problem, the transportation problem,$\dots$) the optimal solution of the problem is invariant to a row and a column permutation of the cost matrix. It is therefore a good question to ask if there is a pair of permutations such that the permuted matrix is Monge. We introduce a generalization of a permutation algorithm by Deineko \& Filonenko~\cite{DeiF79} for real matrices. In O($m^2 + n^2 + mn$) where $m,n$ are the dimensions of the matrix the algorithm decides if there are permutations of rows and columns such that the resulting matrix is from $\ISM$. The question for matrices with weak Monge property is still an open problem and does not seem to have a straightforward correspondence with the algorithm given by Deineko \& Filonenko.
\subsection{Lemmata for the derivation of the algorithm}
	In this section we prove lemmata that are necessary for the derivation of the permutation algorithm. We denote by $\M(\sigma,\pi)$ a matrix $\M$ permuted by a row permutation $\sigma$ and a column permutation $\pi$. If there are permutations $\sigma, \pi$ such that $\M(\sigma,\pi) \in \ISM$ we say that $\M$ is \textit{Monge permutable}.\\
	The first lemma shows that $\ISM$ is closed under an operation of \textit{flipping the matrix upside down and left to right}.
	\begin{Lem}\label{lem6:1}
		Let $\M \in \ISM^{m \times n}$. Define $\sigma(i)=m-i+1$ and $\pi(j)=n-j+1$. Then $\M(\sigma,\pi) \in \ISM$.
	\end{Lem}
	\proof
		For every pair of indices $i,j$ we have that 
		\[
		\ol m_{\sigma(i),\pi(j)} + \ol m_{\sigma(i+1),\pi(i+1)} = \ol m_{m-i+1,n-j+1} + \ol m_{m-i,n-j}.\]
		From Monge property we have
		\[\ol m_{m-i+1,n-j+1} + \ol m_{m-i,n-j} \leq \ul m_{m-i,n-j+1} + \ul m_{m-i+1,n-j},\]
		but the righthand side of the inequality is equal to
		\[
		\ul m_{m-i,n-j+1} + \ul m_{m-i+1,n-j} = \ul m_{\sigma(i),\pi(j+1)} + \ul m_{\sigma(i+1),\pi(j)}.
		\]
		By Theorem~\ref{thm4}.2 we conclude that $\M(\sigma,\pi) \in \ISM$.
	\endproof	
	The following lemma provides a better understanding of what happens if the order of columns is ambiguous meaning that $\ol m_{ij} + \ol m_{k\ell} \leq \ul m_{i\ell} + \ul m_{kj}  \text{ and } \ol m_{i\ell} + \ol m_{kj} \leq \ul m_{ij} + \ul m_{k\ell}$. If this happens, the order of columns and rows does not really matter because all four interval entries are actually real values and so are all entries vertically and horizontally in between them.
	\begin{Lem}\label{lem6:2}
		Let $\M \in \ISM$ and let row indices $i < k$ and column indices $j < \ell$. If it holds that
		\[
		\ol m_{ij} + \ol m_{k\ell} \leq \ul m_{i\ell} + \ul m_{kj}  \text{ and } \ol m_{i\ell} + \ol m_{kj} \leq \ul m_{ij} + \ul m_{k\ell}
		\]
		then for all rows $o$ such that $i \leq o < p \leq k$ it holds
		\begin{enumerate}
			\item $\m_{oj},\m_{pj},\m_{o\ell},\m_{p\ell} \in \mathbb{R}$,
			\item $\m_{oj} + \m_{p\ell} = \m_{o\ell} + \m_{pj}$.
		\end{enumerate}
	\end{Lem}
	\proof
		The following chain of inequalities 
		\[
		\ol m_{ij} + \ol m_{k\ell} \leq \ul m_{i\ell} + \ul m_{kj} \leq \ol m_{i\ell} + \ol m_{kj} \leq \ul m_{ij} + \ul m_{k\ell}
		\]
		turns into a chain of equalities since the first and last members are the same. Taking
		\(
		\ul m_{ij} + \ul m_{k\ell} = \ol m_{ij} + \ol m_{k\ell}
		\)
		and subtracting $\ul m_{ij}$ and $\ol m_{k\ell}$ we have
		\[
		-2 \cdot  m_{k\ell}^{\Delta} = \ul m_{k\ell} - \ol m_{k\ell} = \ol m_{ij} - \ul m_{ij} = 2 \cdot m_{ij}^{\Delta}, 
		\]
		from which
		\(
		- m_{k\ell}^{\Delta} =  m_{i\ell}^{\Delta}.
		\)
		But this means that 
		\(
		\m_{i\ell},\m_{k\ell} \in \mathbb{R}.
		\)
		Similarly, $\m_{i\ell},\m_{kj} \in \mathbb{R}$. This leads to $\m_{ij} + \m_{k\ell} = \m_{i\ell} + \m_{kj}$. \\
		Let now rows $o,p$ be in between rows $i$ and $k$ i.e. $i \leq o < p \leq k$. Since $\M \in \ISM$ the following chain of inequalities holds
		\[
		\m_{ij} + \ol m_{o\ell} \leq \ul m_{oj} + \m_{i\ell} \leq \ol m_{oj} + \m_{i\ell} \leq \m_{k\ell} + \ul m_{oj}.
		\]
		We can rearrange the inequalites in the following way
		\[
		\m_{ij} - \m_{i\ell} \leq \ul m_{oj} - \ol m_{o\ell} \leq 
		\ol m_{oj} - \ul m_{o\ell}  \leq \m_{kj} - \m_{k\ell}
		\]
		and because the first and the last expression equals,
		\[
		\ul m_{oj} - \ol m_{o\ell} = \ol m_{oj} - \ul m_{o\ell}.
		\]
		This leaves us with $-2 \cdot  m_{oj}^{\Delta} = 2 \cdot  m_{o\ell}^{\Delta}$ and from a similar argument as above we have that $\m_{oj},\m_{o\ell} \in \Rz$ and consequently, $\m_{ij} + \m_{o\ell} = \m_{oj} + \m_{i\ell}$. Similarly, we arrive to $\m_{oj} + \m_{p\ell} = \m_{pj} + \m_{o\ell}$.
	\endproof
	The permutation of rows is based upon a combination of conditions which have to be satisfied for the matrix to be strongly Monge. The conditions are taken in a form
	\[
	\ol m_{ij} - \ul m_{i\ell} \leq \ul m_{kj} - \ol m_{k\ell} \text{ where } i < k, j < \ell.
	\]
	For two rows $i,k$ we take into account the first $b$ and the last $B$ columns.
	\begin{Lem}\label{lem6:3}
		Let $\M \in \IR^{m \times n}$. If $\M \in \ISM$ then for every pair of rows $i$ and $k$ such that $i<k$ it holds
		\[
		B\cdot\left(\sum_{j=1}^{b} \ol m_{ij}\right) - b\cdot\left(\sum_{\ell=n-B+1}^{n} \ul m_{i\ell}\right) \leq B\cdot\left(\sum_{j=1}^{b} \ul m_{kj}\right) - b\cdot\left(\sum_{\ell=n-B+1}^{n} \ol m_{k\ell}\right)
		\]
		where $1 \leq b < n-B+1 \leq n$.	
	\end{Lem}
	\proof
		For $i < k$, it holds for every $j$ such that $1 \leq j \leq b$ and every $\ell$ such that $n-B+1 \leq \ell \leq n$ that
		\[
		\ol m_{ij} -  \ul m_{i\ell} \leq \ul m_{kj} - \ol m_{k\ell}.
		\]
		By picking such an inequality for every pair $(j,\ell)\text{ where } j \in \left\{1,\dots,b\right\} \text{ and } \ell \in \left\{n-B+1,\dots,n\right\}$ and adding all these inequalities together we get the formula above.
	\endproof
	The following lemma gives an algorithm to compute the permutations of rows and columns. 
	\begin{Lem}\label{lem6:4}
		Let $\mathbf{u},\mathbf{v} \in \IR^{n}$. Let $\sigma$ be a permutation of $\left\{1,\dots,n\right\} $ such that whenever
		\[ \sigma(i) < \sigma(j) \text{ then } \ol u_i - \ul v_i \leq \ul u_j - \ol v_j.
		\]
		Then in $O(n^2)$ we can compute $\sigma$ or answer that there is no such permutation.
	\end{Lem}
	\proof
		We construct a directed graph $G = (\left\{1,\dots,n\right\},E)$ where $(i,j) \in E$ if $\ol u_i - \ul v_i \leq \ul u_j - \ol v_j$. If there is a pair of vertices $i,j \in G$ without an edge between them, it means that 
		\[\ol u_i - \ul v_i > \ul u_j - \ol v_j \text{ and } \ol u_j - \ul v_j > \ul u_i - \ol v_i\] and by the definition of $\sigma$ no mutual order of these indices yields the permutation; so we stop. From now on, let us suppose that there is at least one edge between all pairs of vertices in graph $G$.\\
		Now let $c_1,\dots,c_k$ be strongly connected components of $G$ such that $t(c_1)<\dots<t(c_k)$, where $t$ is some topological ordering of strongly connected components of $G$.\\
		Now define $\sigma$ as follows. While $\sigma$ is not defined for all indices $i \in \left\{1,\dots,n\right\}$, pick between indices with unspecified $\sigma(i)$ the one for which the topological number of the corresponding strongly connected components containing the vertex $i$ is minimal. Set $\sigma(i)$ as the smallest number from $\left\{1,\dots,n\right\}$ not assigned yet.\\
		To prove that the construction is correct, let $i,j$ be indices such that $\sigma(i) < \sigma(j)$.
		Then either vertices $i,j$ are from the same component or $i$ is from a component with a smaller topological number than the component containing $j$.
		If $i$ and $j$ are from the same component of $G$, it means by the contruction of $G$ that there are edges $(i,j)$ and $(j,i)$ therefore it holds that $\ol u_i - \ul v_i \leq \ul u_j - \ol v_j$.
		If $i$ is in a component with a smaller topological number than $j$, it means that there is an edge $(i,j)$. But the edge $(i,j)$ corresponds to the inequality $\ol u_i - \ul v_i \leq \ul u_j - \ol v_j$.\\
		There exists an algorithm for finding a topological ordering of strongly connected components of a directed graph running in $O(n + m)$ where $n$ equals the number of vertices and $m$ equals the number of edges (see~\cite{Tarjan76}). Since the number of edges $m$ is in the worst case approximately $m \approx n^2$, the algorithm runs in $O(n^2)$. Defining $\sigma$ from the topological ordering $t$ takes $O(n)$, therefore the whole construction takes $O(n^2)$.
	\endproof
		Finally, we prove a lemma about the first step of our algorithm. In this step a pair of rows is determined. The first permutation $\rho$ of the general algorithm is based upon conditions between these two rows. We demand at least two columns to be strictly ordered, otherwise the permutation $\rho$ will have no effect (we want it to prepermute the matrix). Two columns with a strict order are part of two different so called \textit{ambiguity sets}. According to the logical structure we state the lemma in this subsection, however, the notion of \textit{ambiguity sets} necessary in the lemma becomes clear further in the text. We recommend to the reader to first go through the derivation of the algorithm.
	\begin{Lem}\label{lem6:5}
		Let $\M \in \IR^{m \times n}$. Then a problem to decide if there is a row $r$ such that there are two ambiguity sets of columns for rows $1$ and $r$ can be computed in $O(mn)$. If for every row $r$ there is only one ambiguity set of columns, then the matrix has the strong Monge property.
	\end{Lem}
	\proof
		For every row $k$ and for all neighbouring pairs of columns (i.e. $j,j+1$ for $1 \leq j \leq n-1$)  we check if it holds that
		\begin{equation}\label{permeq1}
		\ol m_{1j} - \ul m_{kj} < \ul m_{1,j+1} - \ol m_{k,j+1} \text{ or } \ol m_{1,j+1} - \ul m_{k,j+1} < \ul m_{1j} - \ol m_{kj}.
		\end{equation}
		Only one of these inequlities can hold at the same time because otherwise 
		\[
		\ol m_{1j} - \ul m_{kj} < \ul m_{1,j+1} - \ol m_{k,j+1} \leq \ol m_{1,j+1} - \ul m_{k,j+1} < \ul m_{1j} - \ol m_{kj} \leq \ol m_{1j} - \ul m_{kj}
		\]
		which leads to a contradiction $\ol m_{1j} - \ul m_{kj} < \ol m_{1j} - \ul m_{kj}$. If one of the inequalities holds and the other is $=$, then w.l.o.g. consider
		\[
		\ol m_{1j} - \ul m_{kj} < \ul m_{1,j+1} - \ol m_{k,j+1} \text{ and } \ol m_{1,j+1} - \ul m_{k,j+1} = \ul m_{1j} - \ol m_{kj}.
		\]
		From these two inequalities we can derive that
		\[
		\ol m_{1j} - \ul m_{kj} < \ul m_{1,j+1} - \ol m_{k,j+1} \leq  \ol m_{1,j+1} - \ul m_{k,j+1} = \ul m_{1j} - \ol m_{kj}
		\]
		and therefore $\ol m_{1j} - \ul m_{kj} < \ul m_{1j} - \ol m_{kj}$ which is again a contradiction.\\
		This means that if one inequality holds with $<$ the other must hold with $>$, therefore the order of the columns is strict and they cannot be switched. A strict order of two columns means that these columns cannot be in one ambiguity set, therefore we return row $k$.\\
		It might happen that for every pair $j,j+1$ and for row $k$ neither of the inequalities from \ref{permeq1} is strict. It means that
		\begin{equation}\label{permeq2}
		\ol m_{1j} - \ul m_{kj} \geq \ul m_{1,j+1} - \ol m_{k,j+1} \text{ and } \ol m_{1,j+1} - \ul m_{k,j+1} \geq \ul m_{1j} - \ol m_{kj}.
		\end{equation}
		If both of the inequalities are strict for at least one pair $j,j+1$, it means that no order of columns $j,j+1$ satisfy Monge property and in that case we stop.\\
		If both of the inequalities hold with equality $=$ for all pairs of columns $j,j+1$ in the row $k$, it means that 
		\[
		\ol m_{1j} - \ul m_{kj} = \ul m_{1,j+1} - \ol m_{k,j+1} \leq \ol m_{1,j+1} - \ul m_{k,j+1} = \ul m_{1j} - \ol m_{kj} \leq \ol m_{1j} - \ul m_{kj},
		\]
		therefore $\m_{1j} \m_{kj}$ and also $\m_{k,j+1} \m_{1,j+1}$ are real values and therefore $\m_{1j} - \m_{kj} = \m_{1,j+1} - \m_{k,j+1}$. If this happens for all rows $k$ then the matrix is already Monge because every condition holds with equality.\\
		The last case which remains is when one of the inequalities from \ref{permeq2} is strict $>$ and the second one is equal $=$ for at least one row $r$. Then the order is strict again, because there is only one way to permute these two columns in order to satisfy Monge property. Therefore we return row $r$.
		
		Applying this procedure to each of $m-1$ rows the number of conditions to check is at most $2(n-1)$ for each row. We conclude that the problem can be computed in $O(mn)$.
	\endproof
	
\subsection{Special case algorithm}
We first derive an algorithm for special case interval matrices with nontrivial intervals (i.e. the width of interval is larger than $0$).\\
The algorithm chooses two random rows $i,k$, and according to conditions 
\[ \ol m_{ij} + \ol m_{k\ell} \leq \ul m_{i\ell} + \ul m_{kj}\]
it chooses permutation $\rho$ such that $\rho(j) < \rho(\ell)$ if
\[ \ol m_{ij} - \ul m_{kj} \leq \ul m_{i\ell} - \ol m_{k\ell}.\]
Notice that the permutation $\rho$ is unique. Otherwise both inequalities
\[\ol m_{ij} + \ul m_{kj} \leq \ul m_{i\ell} - \ol m_{k\ell} \text{ and } \ol m_{i\ell} + \ul m_{k\ell} \leq \ul m_{ij} - \ol m_{kj} \]
hold and by Lemma~\ref{lem6:2} the intervals are degenerate (i.e. $\m_{ij},\m_{i\ell},\m_{kj},\m_{k\ell} \in \mathbb{R}$).\\
In the same manner we can now choose columns $\rho(1)$ and $\rho(n)$ and define row permutation $\sigma$ such that $\sigma(i) < \sigma(k)$ if
\[ \ol m_{\sigma(i)\rho(1)} - \ul m_{\sigma(i)\rho(1)} \leq  \ul m_{\sigma(k)\rho(n)} - \ol m_{\sigma(k)\rho(n)}.\]
By Lemma~\ref{lem6:2} the permutation $\sigma$ is unique again.\\
Now if the permuted matrix $\M(\sigma,\rho) \in \ISM$, the algorithm returns $(\sigma,\rho)$. If the matrix does not have the strong Monge property, it means that there are four entries which do not satisfy the corresponding inequality. But since the permutations $\sigma,\rho$ are unique, there is no other permutation of $\M$ that satisfies all necessary conditions.\\
Notice that for special case matrices there are actually two ways to permute the initial matrix $\M$. The first one is pair $(\sigma,\rho)$ and the other is given by Lemma~\ref{lem6:1}.

\subsubsection{Pseudocode of the special case algorithm}
	\begin{Alg}(Special case permutation algorithm)\\
	\begin{algorithm}[H]\label{alg1}
	\KwIn{$\M \in \IR^{mxn}$ an interval matrix with nontrivial intervals}
	\KwOut{''YES'' if $\M$ is Monge permutable together with $\M(\sigma,\rho) \in \ISM$ , ''NO'' otherwise}
	Determine permutation $\rho$ such that
	\[\rho(k) < \rho(\ell) \text{ implies that } \ol m_{1k} - \ul m_{2k} \leq \ul m_{1\ell} - \ol m_{2\ell}.
	\] If no such permutation exists, output ''NO''.\\
	Determine permutation $\sigma$ such that $\sigma(i) < \sigma(k)$ implies that
	\[
	\ol m_{i\rho(1)} - \ul m_{i\rho(n)} \leq \ul m_{k\rho(1)} - \ol m_{k\rho(n)}.
	\]\\
	Check if $\M(\sigma,\rho) \in \ISM$. Output ''YES'' with $\sigma, \rho$ if it does and ''NO'' otherwise.		
	\end{algorithm}
	\end{Alg}
\subsection{General case algorithm}
For general interval matrices the special case algorithm might fail because according to the rule given for the construction of $\sigma$ and $\rho$ the permutations might not be defined unambiguously. Therefore we have to employ a slightly modified algorithm which performs one more permutation.
\subsubsection{Derivation of the algorithm}
Let us suppose that matrix $\M \in \IR^{m\times n}$ is Monge permutable i.e. there are permutations $\sigma$ and $\pi$ such that $\M(\sigma,\pi) \in \ISM$. At first, let us suppose that we already know the permutation $\sigma$ and we would like to derive the permutation $\pi$. We could find for every pair of rows all possible permutations of columns such that Monge property is satisfied for the selected pair and after that choose one permutation that satisfies Monge property for all pairs of rows at the same time.\\
Since this approach is ineffective we construct a permutation only for the pair of the first and the last row. If the permutation cannot be constructed, the matrix is not Monge permutable which is a contradiction with the assumption. Therefore the permutation is either uniquely determined or the order of at least two columns is ambiguous.\\
If the permutation is uniquely determined, it must satisfy Monge property for the rest of row pairs, assuming the matrix is Monge permutable.\\
If the permutation is ambiguous, it means that there are two columns $j,\ell$ such that 
\[
\ol m_{1j} - \ul m_{mj} \leq \ul m_{1\ell} - \ol m_{m\ell} \text{ and } 
\ol m_{1\ell} - \ul m_{m\ell} \leq \ul m_{1j} - \ol m_{mj}.
\]
According to Lemma~\ref{lem6:2}, $\m_{1j},\m_{mj},\m_{1\ell},m_{m\ell} \in \mathbb{R}$ and $\m_{1j} + \m_{m\ell} = \m_{mj} + \m_{1\ell}$ and the conditions also hold if we substitute $1,j$ for any other pair of rows. Again, the permutation satisfies the strong Monge property, otherwise the matrix $\M$ is not Monge permutable.\\
The question that remains is how to determine the permutation $\sigma$. If we knew what the first and the last column in permutation $\pi$ was, we could apply the same idea as for permuting the columns. Therefore, we prepermute the columns by permutation~$\rho$. We choose almost random pair of rows $i,k$ (we further show how) and apply the same rule for permutation as in the special case algorithm i.e.
\[
\rho(j) < \rho(\ell) \text{ implies } \ol m_{ij} + \ol m_{k\ell} \leq \ul m_{i\ell} + \ul m_{kj}.
\]
Because the prepermutation is in general ambiguous, it does not give us the first and the last column. It divides the columns into so called \textit{ambiguity sets}. Two columns are in one ambiguity set if their order cannot be unambiguously determined. Even though the order of columns cannot be determined inside one ambiguity set, for two columns from two different sets the order is strictly given. Therefore, the first and the last ambiguity sets contain the candidates for the first and the last column.\\
Even though we cannot exactly determine the first and the last column we can use a combination of conditions for all columns from the first and the last ambiguity set and base the construction of permutation $\sigma$ upon this combination. Lemma~\ref{lem6:3} provides the condition, i.e. $\sigma(i) < \sigma(k)$ implies
\[
B\cdot\left(\sum_{j=1}^{b} \ol m_{ij}\right) - b\cdot\left(\sum_{\ell=n-B+1}^{n} \ul m_{i\ell}\right) \leq B\cdot\left(\sum_{j=1}^{b} \ul m_{kj}\right) - b\cdot\left(\sum_{\ell=n-B+1}^{n} \ol m_{k\ell}\right)\text{.}
\]
Now the trick is that this process yields an equal permutation to the one based only on the first and the last column.\\
To see this let $\tau$ be the permutation of rows given by the first and the last column in the matrix $\M(\sigma,\pi)$. We want to prove that $\sigma$ is equal to $\tau$.\\
For a contradiction let us suppose that there are two rows $i,k$ such that the order under permutation $\tau(k) < \tau(i)$ and $\sigma(i) < \sigma(k)$ differs. This means that $\ol m_{k1} + \ol m_{in} \leq \ul m_{i1} + \ul m_{in}$ and
\begin{equation}\label{eq1}
    B\cdot\left(\sum_{j=1}^{b} \ol m_{ij}\right) - b\cdot\left(\sum_{\ell=n-B+1}^{n} \ul m_{i\ell}\right) \leq
	B\cdot\left(\sum_{j=1}^{b} \ul m_{kj}\right) - b\cdot\left(\sum_{\ell=n-B+1}^{n} \ol m_{k\ell}\right)\text{.}
\end{equation}
There are two cases to consider. In the first case the order of $k$ and $i$ is unambigous for the permutation $\tau$. Then the inequality is strict $\ol m_{k1} + \ol m_{in} < \ul m_{i1} + \ul m_{in}$. It must also hold for each column $j,\ell$ from the first resp. the last ambiguity set of $\sigma$ that $\ol m_{kj} + \ol m_{i\ell} \leq \ul m_{k\ell} + \ul m_{ij}$ otherwise the matrix is not Monge permutable. But combining all conditions together in the same way as in Lemma~\ref{lem6:3} we achieve a strict inequality
\begin{equation}\label{eq2}
B\cdot\left(\sum_{j=1}^{b} \ol m_{kj}\right) - b\cdot\left(\sum_{\ell=n-B+1}^{n} \ul m_{k\ell}\right) <
B\cdot\left(\sum_{j=1}^{b} \ul m_{ij}\right) - b\cdot\left(\sum_{\ell=n-B+1}^{n} \ol m_{i\ell}\right)\text{.}
\end{equation}
Now the righthand side of~(\ref{eq1}) is less or equal to the lefthand side of~(\ref{eq2}) and the righthand side of~(\ref{eq2}) is less or equal to the lefthand side of~(\ref{eq1}) leading into a contradiction
\[
B\cdot\left(\sum_{j=1}^{b} \ol m_{ij}\right) - b\cdot\left(\sum_{\ell=n-B+1}^{n} \ul m_{i\ell}\right) <
B\cdot\left(\sum_{j=1}^{b} \ol m_{ij}\right) - b\cdot\left(\sum_{\ell=n-B+1}^{n} \ul m_{i\ell}\right)\text{.}
\]
In the second case the order of rows $i$ and $k$ is ambiguous under $\tau$ but this means that switching them does not violate any condition as was discussed before. Therefore even though the permutations $\tau$ and $\sigma$ does not have to be identical, they are equal in the sense that we can use both of them for constructing $\pi$.\\
Last thing to discuss is the construction of prepermutation $\rho$. It is essential for the construction of $\sigma$ to have different candidates for both the first and the last column otherwise the construction fails to determine the order of rows. We need to find a pair of rows $i,j$ such that the permutation divides the columns into at least two ambiguity sets. Lemma~\ref{lem6:5} gives us a way to find these rows.

	\subsubsection{Pseudocode of the general case algorithm}
		\begin{Alg}General case permutation algorithm\\
		\begin{algorithm}[H]\label{alg2}
			\KwIn{$\M \in \IR^{mxn}$}
			\KwOut{''YES'' if $\M$ is Monge permutable together with $\M(\sigma,\pi) \in \ISM$ , ''NO'' otherwise}
			Find a row $r$ such that there are at least two column ambiguity sets for rows $1,r$. If every row has one ambiguity set with row $1$ output ''YES'' with $\sigma = id$ and $\pi = id$. If there is a pair of columns $j,j+1$ which cannot be permuted output ''NO''.\\ 
			Determine permutation $\rho$ such that
			\[\rho(k) < \rho(\ell) \text{ implies that } \ol m_{1k} - \ul m_{jk} \leq \ul m_{1\ell} - \ol m_{j\ell}.
			\] If no such permutation exists, output ''NO''.\\
			Determine $b,B \in \left\{1,..n\right\}$ such that $b$ equals to the size of the first \textit{ambiguity set} of $\rho$ and $B$ equals to the size of the last \textit{ambiguity set} of $\rho$.\\
			Determine row permutation $\sigma\text{ such that } \sigma(i) < \sigma(k)$ implies that
			\[
			B\cdot\left(\sum_{j=1}^{b} \ol m_{ij}\right) - b\cdot\left(\sum_{\ell=n-B+1}^{n} \ul m_{i\ell}\right) \leq B\cdot\left(\sum_{j=1}^{b} \ul m_{kj}\right) - b\cdot\left(\sum_{\ell=n-B+1}^{n} \ol m_{k\ell}\right).
			\]
			If no such permutation exists, output ''NO''.\\
			Determine column permutation $\pi$ such that \[\pi(k) < \pi(\ell) \text{ implies that } \ol m_{\sigma(1),k} - \ul m_{\sigma(n),k} \leq \ul m_{\sigma(1),\ell} - \ol m_{\sigma(n),\ell}.\]If no such permutation exists, output ''NO''.\\
			Check if $\M(\sigma,\pi) \in \ISM$. Output ''YES'' with $\sigma, \pi$ if it does and ''NO'' otherwise.		
		\end{algorithm}
	\end{Alg}
\subsection{Complexity of the algorithm}
	The correctness of both variants of the algorithm follows from the derivations. It remains to determine the time complexity of the algorithm.
	\begin{The}\label{thm05:1}
		For $\M \in \mathbb{IR}^{m\times n}$, Algorithm~\ref{alg1} runs in $O(m^2 + n^2 + mn)$.
	\end{The}
	\proof
		By Lemma~\ref{lem6:4} the permutation $\rho$ can be constructed in $O(n^2)$ and $\sigma$ in $O(m^2)$. Using Theorem~\ref{thm4}.(3) it can be checked in $O(mn)$ time if the permuted matrix is strongly Monge. Altogether, the time complexity is $O(m^2 + n^2 + mn)$.
	\endproof
	
	\begin{The}\label{thm05:2}
		For $\M \in \mathbb{IR}^{m\times n}$, Algorithm~\ref{alg2} runs in $O(m^2 + n^2 + mn)$.
	\end{The}
	\proof
	The determination of row $r$ takes $O(mn)$ time according to Lemma~\ref{lem6:5}.
	Construction of both column permutation $\rho,\pi$ takes $O(n^2)$ time and the construction of $\sigma$ takes $O(m^2)$ time as can be seen from a slight modification of Lemma~\ref{lem6:4}. We can easily derive $b$ and $B$ from $\rho$ by checking mostly $2n$ conditions, therefore the determination of $b,B$ runs in $O(n)$.
	Finally, by Theorem~\ref{thm4}.(3) Monge recognition procedure takes $O(mn)$.
	Altogether, the time complexity of Algorithm~\ref{alg2} is $O(m^2 + n^2 + mn)$.
	\endproof
\section{Conclusion}
We introduced two classes of interval Monge matrices - $\ISM$ and $\IWM$. For $\ISM$, following mostly results of real Monge matrices, we generalized several characterizations. For $\IWM$ we offered a polynomial characterization and several necessary and sufficient conditions. In Theorem~\ref{thm11} we indicated a larger class of conditions that might be interesting to further investigate.\\
We presented lists of closure properties under operations on $\ISM$ and $\IWM$ and under operations combining both classes of matrices.\\
We introduced generalization of Deineko \& Filonenko permutation algorithm for interval matrices, which determines if the interval matrix is Monge permutable, i.e. there is a permutation of rows and columns such that the permuted matrix has the strong Monge property.

{\small
{\em Authors' address}:
{\em Martin \v{C}ern\'y}, Charles University, Prague, Czech Republic
 e-mail: \texttt{cerny@\allowbreak kam.mff.cuni.cz}.
}

\end{document}